\documentclass[a4paper,notitlepage,11pt]{article}
\usepackage{amsfonts}
\usepackage{amsmath}
\usepackage{amssymb}
\usepackage{amsopn}
\usepackage[latin1]{inputenc}
\usepackage{graphicx}
\usepackage[lmargin=2.5cm, rmargin=2.5cm,bottom=2.5cm,top=2.5cm]{geometry}
\newtheorem{theorem}{Theorem}
\newtheorem{lemma}[theorem]{Lemma}

\newtheorem{example}{Example}
\newtheorem{proposition}[theorem]{Proposition}
\newtheorem{question}{Question}
\newtheorem{remark}{Remark}

\newcommand{\nats}{{\mathbb N}}

\title {On minimal factorizations of words as products of palindromes}

\author
{
    A. Frid\thanks{Corresponding author, Sobolev Institute of Mathematics SB RAS;
supported in part by Presidential grant MK-4075.2012.1 and by RFBR grant 12-01-00089.},
    S. Puzynina \thanks{Sobolev Institute of Mathematics SB RAS and University of Turku; supported in part by Academy of Finland (grant 251371) and by RFBR (grants 10-01-00424 and 12-01-00448)}, L. Zamboni \thanks{Universit\'e Lyon I and University of Turku, supported in part by a FiDiPro grant from the Academy of Finland.}
}

\begin{document}
\maketitle
\begin{abstract}
 Given a finite word $u$, we define its \emph{palindromic length} $|u|_{\rm pal}$ to be the least number $n$ such that $u=v_1v_2\dots v_n$ with each $v_i$ a palindrome.
 We address the following open question: Does there exist an infinite non ultimately periodic word $w$ and a positive integer $P$ such that $|u|_{\rm pal} \leq P$
 for each factor $u$ of $w$? In this note we give a partial answer to this question. Let $k$ be a positive integer. We prove that if an infinite word $w$ is $k$-power free, then for each positive integer $P$ there exists a factor $u$ of $w$ whose palindromic length $|u|_{\rm pal}>P.$ We also extend this result to a wider class of words satisfying the so-called $(k,l)$-condition, which includes for example the Sierpinski word.
\end{abstract}

\section{Introduction}
Let $A$ be a finite non-empty set, and let $A^+$ denote the set of all finite non-empty words in $A.$ A word $u=u_1u_2\cdots
u_n\in A^+$ is called a {\it palindrome}  if $u_{i}=u_{n-i+1}$
for each $i=1,\ldots,n-1$.  In particular each $a\in A$ is a palindrome. We also regard the empty word as a palindrome. 

Palindrome factors of finite or infinite words have been studied from different points of view. In particular, Droubay, Justin and Pirillo \cite{djp} proved that a word of length $n$ can contain at most $n+1$ distinct palindromes, which gave rise to the theory of {\it rich} words. The number of palindromes of a given length occurring in an infinite word is called its  {\it palindrome complexity} and is bounded by a function of its usual subword complexity \cite{abcd}. However, in this paper we study palindromes in an infinite word from the point of view of decompositions.

 For each  word $u\in A^+$ we define its \emph{palindromic length}, denoted by $|u|_{\rm pal}$, to be the least number $P$ such that $u=v_1v_2\cdots v_P$ with each $v_i$ a palindrome. As each letter is a palindrome, we have $|u|_{\rm pal}\leq |u|$, where $|u|$ denotes the length of $u.$
For example, $|01001010010|_{\rm pal}=1$ while  $|010011|_{\rm pal}=3.$ Note that $010011$ may be expressed as a product of $3$ palindromes in two different ways: $(0)(1001)(1)$ and $(010)(0)(11).$
In \cite{R}, O.~Ravsky obtains an intriguing formula for the supremum of the palindromic lengths of all binary words of length $n.$

The question considered in this paper is
\begin{question}\label{c1}
Does there exist an infinite non ultimately periodic word $w$ and a positive integer $P$ such that $|u|_{\rm pal} \leq P$
 for each factor $u$ of $w$?
\end{question}

We conjecture that such aperiodic words do not exist, but at the moment we can prove it only partially. Namely, in this paper we prove that if such a word exists, then it is not $k$-power-free for any $k$ and moreover, for all $k>1$, $l\geq 0$ it does not satisfy the $(k,l)$-condition defined in Section \ref{s:kl}. A discussion what exactly the condition means and which class of words should be studied now to give a complete answer to the question is given in Section \ref{s:disc}. 

A preliminary version of this paper has been reported at Journ\'ees Montoises 2012.

\section{The case of $k$-power-free words}
Let $k$ be a positive integer. A word $v\in A^+$ is called a  {\it $k$-power} if $v=u^k$ for some word $u\in A^+$. An infinite word $w=w_1w_2\ldots \in A^{\nats}$ is said to be {\it $k$-power-free} if no factor $u$ of $w$ is a $k$-power. For instance, the Thue-Morse word $0110100110010110\ldots$ fixed by the morphism $0\mapsto 01,$ $1\mapsto 10$ is $3$-power free \cite{lothaire}.

\begin{theorem}\label{t:kpf} Let $k$ be a positive integer and $w=w_1w_2\ldots \in A^\nats.$ If $w$ is $k$-power free, then for each positive integer $P$ there exist a prefix $u$ of $w$ with
$|u|_{\rm pal}>P$.
\end{theorem}

Recall that a word $u_1\cdots u_n$ is called $t$-periodic if $u_i=u_{i+t}$ for all $1\leq i \leq n-t$.

The proof of Theorem~\ref{t:kpf} will make use of the following lemmas.

\begin{lemma}\label{l:period}
Let $u$ be a palindrome. Then for every palindromic proper prefix $v$ of $u$, we have that $u$ is $(|u|-|v|)$-periodic.
\end{lemma}

{\sc Proof.}  If $u$ and $v$ are palindromes with $v$ a proper prefix of $u$, then $v$ is also a suffix of $u$ and hence $u$ is  $(|u|-|v|)$-periodic.
\hfill $\Box$

In what follows, the notation $w[i..j]$ can mean the factor $w_i w_{i+1}\cdots w_j$ of a word $w=w_1\cdots w_n \cdots$ as well as its precise occurrence starting with the position numbered $i$; we always specify it when necessary.

\begin{lemma}\label{l:1+1/k}
Suppose the infinite word $w$ is $k$-power-free. If $w[i_1..i_2]$ and $w[i_1..i_3]$ are palindromes with $i_3>i_2$, then
\[\frac{|w[i_1..i_3]|}{|w[i_1..i_2]|}> 1+\frac{1}{k-1}. \hfill \Box \]
\end{lemma}
{\sc Proof.} By Lemma \ref{l:period}, the word $w[i_1..i_3]$ is $(i_3-i_2)$-periodic; at the same time, it cannot contain a $k$-power, so, $|w[i_1..i_3]|<k(i_3-i_2).$
Thus,
\[\frac{|w[i_1..i_3]|}{|w[i_1..i_2]|}=\frac{|w[i_1..i_3]|}{|w[i_1..i_3]|-(i_3-i_2)}>
 \frac{|w[i_1..i_3]|}{\left (1-\frac{1}{k}\right )(|w[i_1..i_3]|)}=1+\frac{1}{k-1}. \hfill \Box
\]

\begin{lemma}
Let $N$ be a positive integer. Then for each $i\geq 0,$  the number of palindromes of the form $w[i..j]$  of length less or equal to $N$   is at most $2+\log_{k/(k-1)}N$.
\end{lemma}
{\sc Proof.} For each $i\geq 0,$ the length of the shortest non-empty palindrome beginning in position $i$ is equal to $1.$ By the previous lemma, the next palindrome beginning in position $i$ is of length more than $\frac{k}{k-1}$, and the one after that is of length more than $(\frac{k}{k-1})^2$, and so on.  The longest one is of length at most $N$ but more than
$(\frac{k}{k-1})^P$, so that $P\leq \log_{k/(k-1)}N$, and the total number $n+1$ of such words is at most $1+\log_{k/(k-1)}N.$ Adding the empty word which is a palindrome gives the desired result.\hfill $\Box$

\smallskip
{\sc Proof of Theorem \ref{t:kpf}.}  Fix a positive integer $P$ and let $N$ be a positive integer satisfying\[(2+\log_{k/(k-1)}N)^P<N.\]
By the previous lemma, the number of prefixes of $w$ of the form $v_1v_2\ldots v_P$, where each $v_i$ is a palindrome, of length less or equal to $N$ is at most $(2+\log_{k/(k-1)}N)^P,$ and hence at most $N.$   But $w$ has $N$-many non-empty prefixes of length less or equal to $N.$ This means that there exists a prefix $u$ of $w$ of length less or equal to $N$ such that
$|u|_{\rm pal}>P$.
\hfill $\Box$

\section{Privileged words and other regularities}\label{s:priv}
In fact, the proof above does not use directly any properties of palindromes except for Lemma~\ref{l:1+1/k}. So, analogous statements on the properties of $k$-power-free words can be proved for any other type of word regularities for which lemmas analogous to Lemma \ref{l:1+1/k} hold. In particular, we can almost immediately extend Theorem \ref{t:kpf} to {\it privileged words}.

Privileged words have been introduced by J. Kellendonk, D. Lenz and J. Savinien \cite{K}; they are studied also in \cite{pelto}. Privileged words are defined recursively as follows: first, the empty word and each element $a\in A$ is privileged. Next, a word $u\in A^+$ with $|u|\geq 2$ is privileged if and only if it is a complete first return to a shorter non-empty privileged word, i.e., if there exists a non-empty privileged word $v$ which is both a proper prefix and a proper suffix of $u$ and which occurs  in $u$ exactly twice.
For example, $00$ is privileged as it is a complete first return to the privileged word $0.$ Similarly, $00101100$ is privileged as it is a complete first return to the privileged word $00.$ This latter example shows that a privileged word need not be a palindrome. Conversely, the palindromes $1231321$ and $00101100110100$ are not privileged as neither word is a complete first return. However, if a word $w$ is {\it rich} meaning that it contains $|w|+1$ factors which are palindromes, then each non-empty factor of $w$ is a palindrome if and only if it is privileged (see Proposition 2.3 in \cite{K}).  Thus analogously we define the {\it privileged length} of a word $u\in A^+,$ denoted $|u|_{\rm priv}$ to be the least number $n$ such that $u=v_1v_2\cdots v_n$ with each $v_i$ a privileged word.
Again we have the inequality $|u|_{\rm priv}\leq |u|.$
For instance, $|00101100|_{\rm priv}=1,$ while  $|00101100110100|_{\rm priv}=3.$ We note that $00101100110100$ may be written as a product of $3$ privileged words in more than one way: $(0)(010110011010)(0)$ or $(00)(1011001101)(00).$ 

The following lemma is analogous to Lemma \ref{l:period}.
\begin{lemma}\label{l:period_priv}
Let $u$ be a privileged word. Then for every privileged proper prefix $v$ of $u$, we have that $u$ is $(|u|-|v|)$-periodic.
\end{lemma}
{\sc Proof.} Suppose $u$ and $v$ are privileged words with $v$ a proper prefix of $u.$ We will prove that
$v$ is also a suffix of $u.$ We proceed by induction on $|u|.$ The result is vacuously true for $|u|=1.$ Next suppose $|u|>1.$ Then $u$ is a complete first return to a privileged word $u'$ with $|u'|<|u|.$ 
We claim that $|v|\leq |u'|.$ In fact, suppose to the contrary that $|v|>|u'|.$ Then $u'$ would be a proper prefix of $v$ and hence by induction hypothesis $u'$ is also a suffix of $v.$ This means that $u'$ occurs at least three times within $u$ (as a prefix of $v,$ as a suffix of $v$ and as a suffix of $u$). This contradicts that $u$ is a complete first return to $u'.$
Having established that $|v|\leq |u'|,$ it follows that $v$ is a suffix of $u'.$ In fact, if $|v|=|u'|,$ then $v=u'$ while if $|v|< |u'|,$ then by induction hypothesis $v$ is a suffix of $u'.$ As $u'$ is a suffix of $u$ we obtain that $v$ is a suffix of $u$ as required. Whence, $u$ is  $(|u|-|v|)$-periodic.
\hfill $\Box$

Now using Lemma \ref{l:period_priv} we can prove the following lemma.

\begin{lemma}\label{l:1+1/k_priv}
Suppose the infinite word $w$ is $k$-power-free. If $w[i_1..i_2]$ and $w[i_1..i_3]$ are  privileged words with $i_3>i_2$, then
\[\frac{|w[i_1..i_3]|}{|w[i_1..i_2]|}\geq 1+\frac{1}{k-1}.\]
\end{lemma}
{\sc Proof.} By Lemma \ref{l:period_priv}, the word $w[i_1..i_3]$ is $(i_3-i_2)$-periodic; the rest of the proof is completely analogous to that of Lemma \ref{l:1+1/k}. \hfill $\Box$

Now, using this lemma, we analogously to Theorem \ref{t:kpf} prove the following
\begin{theorem}\label{t:kpf_priv} Let $k$ be a positive integer and $w=w_1w_2\ldots \in A^\nats.$ If $w$ is $k$-power free, then for each positive integer $P$ there exist a prefix $u$ of $w$ with
$|u|_{\rm priv}>P$.
\end{theorem}

Instead of privileged words, we could use words of any other type for which a statement analogous to Lemmas \ref{l:1+1/k} and \ref{l:1+1/k_priv} would hold. However, the proof of the next more general statement uses substantially the properties of palindromes and not any other family of words.

\section{The case of the $(k,l)$-condition} \label{s:kl}

Recall that a fractional power $w^{p/q}$ of a word $w$ whose length $|w|$ is divisible by $q$ is defined as the word $w^{\lfloor p/q \rfloor}w'$, where $w'$ is the prefix of $w$ of length $\{p/q\}|w|$. To state the next, more general case for which we can prove the unboundedness of the palindromic length, let us fix some integer $k>0$ and define a \emph{$k$-run} in a word $w$ as follows: a $k$-run is an occurrence $w[i..j]$ such that the word $w[i..j]$ is a $k'$-power for some (possibly fractional) $k'\geq k$, but neither $w[i-1..j]$ nor $w[i..j+1]$ are $k''$-powers for any $k''\geq k$. 

Several problems on the maximal number and the sum of exponents of runs in a finite word have been studied, e.~g., in \cite{cit,ckrrw}.

\begin{example}
 {\rm
The word $v=1(100)^5101(001)^711$ has two $5$-runs, namely, $v[2..18]=(100)^{17/3}$ and $v[18..40]=(010)^{23/3}$. They are also 2-runs, 3-runs and 4-runs; the latter is also a 6-run and a 7-run. There are no 8-runs in $v$.
}
\end{example}

Note that $k$-runs are defined as occurrences of words, or, in fact, as pairs of positions corresponding to their beginnings and ends. So, we can say that a $k$-run $w[i..j]$ covers a position $x$ of the word $w$ if $i+1 \leq x \leq j-1$. Clearly, a position can be covered by an arbitrary number of $k$-runs, and even by an infinite number of them.

\begin{example}
  {\rm
The third letter in the word $v=00010010010100100101001001$ is covered by three 3-runs: $v[1..3]=0^3$, $v[2..11]=(001)^{10/3}$ and $v[3..26]=(01001001)^3$.

The position 1 in the infinite word defined as the limit of the sequence $\{s_i\}_{i=1}^\infty$, where $s_1=0001$, $s_{i+1}=s_i^3 1$, is covered by an infinite number of 3-runs.
}
\end{example}
Let us denote the number of $k$-runs covering position $n$ in a word $w$ by $r_{w,k}(n)$ or simply by $r_k(n)$ if $w$ is fixed in advance. The maximal value of $r_k(n)$ for $i \leq n \leq j$ is denoted by $r_k[i..j]$.

We say that an infinite word $w$ satisfies the {\em $(k,l)$-condition} for some $k\geq 2$ and $l\geq 0$ if it is not ultimately periodic and $r_{w,k}(n)\leq l$ for all $n$, that is, if each position $n$ in $w$ is covered by at most $l$ many $k$-runs.

\begin{example}
 {\rm The Sierpinski word $w_s=0101110101^90101110101^{27}\cdots$, defined as the fixed point starting with 0 of the morphism $\varphi: 0 \mapsto 010, 1 \mapsto 111$,  satisfies the $(3,1)$-condition. Indeed, the only primitive factor $u$ whose powers at least 3 occur in $w_s$ is $1$, and thus there is at most one 3-run covering each position in it.
}
\end{example}

\begin{remark}\label{r1}
{\rm
We have not managed to find a proof of this statement in the literature, but it seems very believable that for any morphic word $w$ there exists some $k$ such that there exists only a finite number of primitive words whose powers greater than $k$ occur in $w$. If it is true, it means almost immediately that all the morphic words satisfy a $(k,1)$-condition for some $k$, and thus that all aperiodic morphic words have unbounded palindromic length of factors.  
}
\end{remark}

The following theorem is a generalization of Theorem \ref{t:kpf} which corresponds to the particular case of $l=0$. Note that it is stated for a factor of the word $w$, not for its prefix.

\begin{theorem}\label{t:klc}
 If an infinite word $w$ satisfies the $(k,l)$-condition for some $k\geq 2$ and $l \geq 0$, then for each given $P>0$ it contains a factor $u$ with $|u|_{\rm pal}>P$.
\end{theorem}

The remaining part of the section is devoted to the proof of this
theorem.  The scheme of the proof is the following. In the first
part of the proof we are going to introduce a new measure of
words, so that the ratio of measures of two palindromes
starting from a point is at least $1+\frac{C}{k-1}$ for some
constant $C$ (Lemma \ref{l:1+c/k}). In the second part of the
proof we choose a factor of big enough measure and deduce that the
factorizations of prefixes of this factor into $P$ palindromes
cannot cover all the prefixes. Hence we derive the existence of a
factor with palindromic length greater than $P$. Though the
general idea of the proof is similar to the case of $k$-power free
words (with measure instead of length), the proof is much more
technical.

\smallskip

First of all, note that if there exist arbitrarily long parts $w[i..j]$ of $w$ with $r_k[i..j]=0$, then we can proceed as in the proof of Theorem \ref{t:kpf} and find a prefix of palindromic length greater than $P$ in any factor $w[i..i+N]$ of $w$ such that $r_k[i..i+N]=0$ and $(2+\log_{k/(k-1)}N)^P<N$. So, in the main case of the theorem the length of factors of $w$ not intersecting with any $k$-runs is bounded, and in particular $w$ contains an infinite number of $k$-runs. So, from now on we assume that $w$ satisfies this condition.

Let us say that a $k$-run $w[i..j]$ is an {\it upper} $k$-run in
$w$ if it is not covered by another $k$-run, that is, if there is
no $k$-run $w[i'..j']$ not equal to it such that $i'\leq i<j\leq
j'$. For each position $n$ in $w$, we define by $m_k(n)=m(n)$ the number
of upper $k$-runs of the form $w[n..i]$ or $w[i..n]$.

\begin{example}
 {\rm
Consider the word $v=11(1000)^3 0$. In it, $v[1..3]=111$, $v[3..14]=(1000)^3$ and $v[12..15]=0000$ are upper 3-runs, whereas the 3-runs $v[4..6]=000$ and $v[8..10]=000$ are not upper since they are covered by $v[3..14]=(1000)^3$. Also, we have $m(3)=2$, $m(1)=m(12)=m(14)=m(15)=1$ and $m(n)=0$ for any other $n$; here $m(n)$ is exactly the number of occurrences of the position $n$ in the notation of the upper 3-runs, which are  $v[1..3]$, $v[3..14]$ and $v[12..15]$.
}
\end{example}

\begin{lemma}\label{l:02}
For each $n$ we have $0 \leq m(n) \leq 2$.
\end{lemma}
{\sc Proof.} At each position $n$ of $w$ there is at most one upper $k$-run beginning in position $n$ and at most one upper $k$-run ending in position $n$, so, $0 \leq m(n) \leq 2$.  \hfill $\Box$

\smallskip
Now let us define the {\it measure} $m[i..j]$ of an occurrence $w[i..j]$ as the sum
\[m[i..j]=\sum_{n=i}^{j}m(n).\]
Clearly, it is indeed a measure, that is, the function $m[i..j]$ is non-negative, equal to 0 for the empty word $w[i..i-1]$ for any $i$, and the measure of a disjoint union is the sum of measures, which means in particular that 
\begin{lemma}\label{l:measure}
 For all $i_1\leq i_2 < i_3$ we have
\[m[i_1..i_3]=m[i_1..i_2]+m[i_2+1..i_3]. ~ \Box\]
\end{lemma}
Note also that the function $r_k[i..j]$ is defined as the maximum of $r_k(n)$ for $n \in \{i,\ldots,j\}$, and is uniformly bounded by $l$ due to the $(k,l)$-condition, whereas $m[i..j]$ is defined as the sum of $m(n)$ and is not uniformly bounded since otherwise $w$ would contain a finite number of $k$-runs. Recall that we assume that the length of factors of $w$ not intersecting with any $k$-runs is uniformly bounded, since otherwise we simply apply the proof of Theorem \ref{t:kpf}. This gives us

\begin{lemma}\label{l:l'M}
 There exists a unique $l'$, $1 \leq l' \leq l$, such that
\begin{itemize}
 \item there exists some $M>0$ such that $m[i..j]\geq M$ implies $r_k[i..j]\geq l'$;
 \item for all $L$ there exist $i,j$ such that $m[i..j]\geq L$ and $r_k[i..j]\leq l'$. \hfill $\Box$
\end{itemize}
\end{lemma}
In the remaining part of the proof we shall always consider (occurrences of) factors of $w$ with $r_k[i..j]\leq l'$: due to the lemma above, their measure can be arbitrarily large. Due to the same lemma, each part of such a word whose measure is at least $M$ must contain a position covered by exactly $l'$ of $k$-runs.

Let us say that a $k$-run $w[j_1..j_2]$ is an {\em internal} run within an occurrence $w[p_1..p_2]$ if it intersects it but does not cover positions $p_1$ or $p_2$, that is, if $p_1 < j_1<j_2 < p_2$. Similarly, a $k$-run $w[j_1..j_2]$ is called a {\em left} $k$-run in a word $w[p_1..p_2]$ if it covers the position $p_1$ but not $p_2$, that is, $j_1 \leq p_1 \leq j_2<p_2$; symmetrically, it is a {\em right} $k$-run if it covers the position $p_2$ but not $p_1$, that is, $p_1<j_1 \leq p_2 \leq j_2$. At last, it is a {\em covering} $k$-run if it covers both ending positions, that is, if $j_1 \leq p_1 < p_2 \leq j_2$.

\begin{lemma}\label{l:runs1}
Each $k$-run intersecting with an occurrence $w[p_1..p_2]$ is either internal, or left, or right, or covering.
\end{lemma}
{\sc Proof.}
The four cases are determined by two facts: if the symbols $w_{p_1}$ and $w_{p_2}$ are parts of the $k$-run. \hfill $\Box$

\begin{lemma}
If $w[p_1..p_2]=w[q_1..q_2]$ and $w[p_1+i..p_2-j]$ for some $i,j>0$ is an internal $k$-run for $w[p_1..p_2]$, then $w[q_1+i..q_2-j]$ is an internal $k$-run in $w[q_1..q_2]$. However, the sets of left, right and covering runs depend on an occurrence of a word, and thus can be completely different for $w[p_1..p_2]$ and for $w[q_1..q_2]$. Moreover, if  $w[p_1+i..p_2-j]$ is an upper $k$-run, it does not imply that $w[q_1+i..q_2-j]$ is an upper $k$-run, and vice versa.
\end{lemma}
{\sc Proof.}If $w[p_1..p_2]=w[q_1..q_2]=u=u_1\cdots u_n$, where $n=p_2-p_1+1=q_2-q_1+1$, then $w[p_1+i..p_2-j]=u_{i+1}\cdots u_{n-j}=v$. By the maximality in the definition of a $k$-run, we see that the symbols $u_i$ and $u_{n-j+1}$ break the periodicity of $v$, so the $k$-run always starts at the symbol number $i+1$ of $u$ and ends at its symbol number $n-j+1$.

To give an example for the second part of the statement, consider the word
\[aba.bb.aba.a^{k-2}.aba.bb.(aba)^{k}cw',\]
where the infinite word $w'$ is on the alphabet $\{b,c\}$. We see that in the first occurrence of $aba$ in it, there are no left, right or covering $k$-runs; in the second one, there is a right $k$-run $a^k$, which is also a left run for the third occurrence of $aba$; and the fourth occurrence of $aba$ is a prefix of a $k$-run $(aba)^k$ covering it.

At last, one occurrence of a word can be an upper internal $k$-run whereas another occurrence is not an upper one. As an example, consider the word
\[(ba^kb)^k bb \; ba^kb c w',\]
where $w'$ does not contain the symbol $a$. We see that in the first $k$ occurrences of $b a^k b$ the $k$-runs $a^k$ are covered by the $k$-run $(ba^k b)^kb$, and in the last one, the $k$-run $a^k$ is an upper one. \hfill $\Box$

\smallskip
To state the next lemma, symmetric to the first part of the previous one, we denote by $\tilde{v}$ the mirror image $v_n v_{n-1} \cdots v_1$ of a word $v=v_1\cdots v_n$. In particular, a palindrome is exactly a word $v$ such that $v=\tilde{v}$.
\begin{lemma}\label{l:vtilde}
If $w[p_1..p_2]=v$ and $w[q_1..q_2]=\tilde{v}$, and if $w[p_1+i..p_2-j]$ for some $i,j>0$ is an internal $k$-run within $w[p_1..p_2]$, then $w[q_1+j..q_2-i]$ is an internal $k$-run in $w[q_1..q_2]$.
\end{lemma}
{\sc Proof.} It is sufficient to realize that if $v[i..j]$ is an internal $k$-run within $v=v_1\cdots v_n$, then $\tilde{v}[n-j+1..n-i+1]$ is an internal $k$-run withing $\tilde{v}=\tilde{v}_1\cdots \tilde{v}_n=v_n\cdots v_1$. \hfill $\Box$

\begin{lemma}
The number of left (resp., right) $k$-runs for $w[p_1..p_2]$ is bounded by $r_(p_1)$ (resp., $r_k(p_2)$) and thus by $r_k[p_1..p_2]$.
\end{lemma}
{\sc Proof.} The $k$-runs covering the position $p_1$ are left or covering runs for $w[p_1..p_2]$, and their total number is $r_k(p_1)$. Symmetrically, the $k$-runs covering the position $p_2$ are right or covering runs for $w[p_1..p_2]$, and their total number is $r_k(p_2)$. \hfill $\Box$

\begin{lemma}\label{l:covering}
The set of covering runs for $w[p_1..p_2]$ can be non-empty only if $m[p_1..p_2]\leq 2r_k[p_1..p_2]$.
\end{lemma}
{\sc Proof.}
Suppose that there is a $k$-run covering $w[p_1..p_2]$. Then every $k$-run $w[q_1..q_2]$ 
contributing to the measure of $w[p_1..p_2]$ satisfy the following condition: one of the values $q_1, q_2$ lies inside
the interval $[p_1..p_2]$, the other one outside. In other words, either
$p_1< q_1$, $p_1\leq q_2\leq p_2$, or $p_1\leq q_1\leq p_2$, $p_2<
q_2$. The number of $k$-runs satisfying the first condition is at
most $r_k(p_1)$, the number of $k$-runs satisfying the second
condition is at most $r_k(p_2)$, and each of them contributes at most $1$ to
the measure. Therefore, $m[p_1..p_2]\leq r_k(p_1)+r_k(p_2)\leq
2r_k[p_1..p_2]$. \hfill $\Box$

\smallskip
The following two ``lemmas of inviolable parts'' are crucial for the proof. The first one is stated for an occurrence of $\tilde{v}$ and the second one for an occurrence of $v$, since it is what we need further in the proof, but in fact both of them can be stated both for $v$ and for $\tilde{v}$.

\begin{lemma}\label{l:ivp1}
If $w[i_1..i_2]=v$, $r_k[i_1..i_2]=l'$, and $m[i_1..i_2]\geq M+4l'$,
where the parameters $l'$  and $M$ are defined in Lemma
\ref{l:l'M}, then $m[j_1..j_2] \geq 2$ for all $j_1, j_2$ such
that $w[j_1..j_2]=\tilde{v}$ and $r_k[j_1..j_2]\leq l'$.
\end{lemma}
{\sc Proof.} Due to Lemma \ref{l:covering}, there are no $k$-runs covering $w[i_1..i_2]$. Consider the maximal $q_1\in \{i_1,\ldots,i_2-1\}$ such that $w[x..q_1]$ is a left $k$-run for $w[i_1..i_2]$; this $k$-run is a covering one for $w[i_1..q_1]$, and thus due to Lemma \ref{l:covering}, $m[i_1..q_1]\leq 2r_k[i_1..q_1]\leq 2r_k[i_1..i_2]=2l'$.  If there are no left $k$-runs in $w[i_1..i_2]$, we put $q_1=i_1-1$, and thus $m[i_1..q_1]=0$ since it is an empty word. So, due to Lemma \ref{l:measure}, we have $m[q_1+1..i_2]=m[i_1..i_2]-m[i_1..q_1]\geq M+2l'$.

Now symmetrically, due to Lemma \ref{l:covering}, there are no $k$-runs covering $w[q_1+1..i_2]$. Consider the minimal $q_2\in \{i_1+1,\ldots,i_2\}$ such that $w[q_2..y]$ is a right $k$-run for $w[i_1..i_2]$. In fact we have $q_2\geq q_1+2$ since otherwise $w[q_2..y]$ is a covering $k$-run for $w[q_1+1..i_2]$, a contradiction. So, the $k$-run $w[q_2..y]$ is a right $k$-run for $w[q_1+1..i_2]$ and  a covering one for $w[q_2..i_2]$, and thus due to Lemma \ref{l:covering}, $m[q_2..i_2]\leq 2r_k[q_2..i_2]\leq 2r_k[i_1..i_2]=2l'$.  If there are no right $k$-runs in $w[i_1..i_2]$, we put $q_2=i_2+1$, and thus $m[q_2..i_2]=0$ since it is an empty word. So, due to Lemma \ref{l:measure}, we have $m[q_1+1..q_2-1]=m[q_1+1..i_2]-m[q_2..i_2]\geq M$.

Consider the occurrence $w[q_1+1..q_2-1]$. By the construction, all $k$-runs covering it are internal $k$-runs for $w[i_1..i_2]$. At the same time, $m[q_1+1..q_2-1]\geq M$, and thus by the definition of $M$ in Lemma \ref{l:l'M}, there is a position
$i_1+n \in \{q_1+1,\ldots,q_2-1\} $ with $r_k(i_1+n)=l'$.
But due to Lemma \ref{l:vtilde}, all the $l'$ internal $k$-runs for $w[i_1..i_2]=v$ covering the position $i_1+n$ have analogues in $w[j_1..j_2]=\tilde{v}$ which cover the position $j_2-n \in \{j_1,\ldots,j_2\}$. Due to the condition $r_k[j_1..j_2]\leq l'$, at least one of these $l'$ runs internal for $w[j_1..j_2]$ and covering the position $j_2-n$ is an upper run in $w$, and thus it contributes 2 to $m[j_1..j_2]$. \hfill $\Box$.

\begin{figure}
\begin{center}
\includegraphics[width=0.8\textwidth]{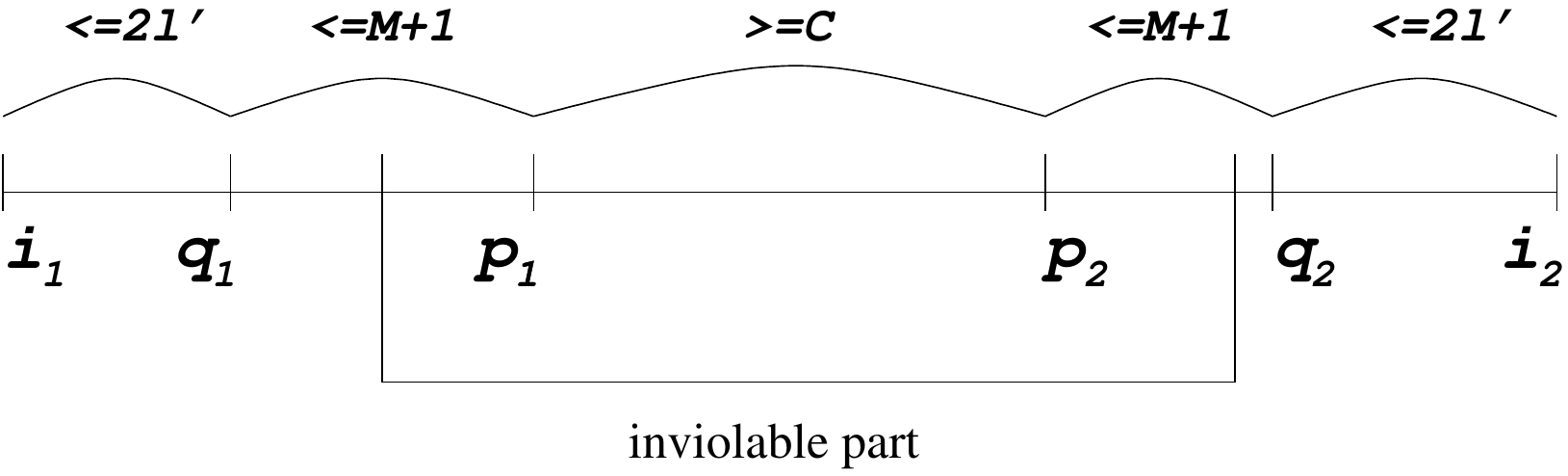}\label{f:f1}
\end{center}
\caption{Proof of Lemma \ref{l:ivp2}}
\end{figure}

\begin{lemma}\label{l:ivp2}
 If $w[i_1..i_2]=v$, $r_k[i_1..i_2]=l'$, and $m[i_1..i_2]= 2M+4l'+2+C$, where the parameters $l'$  and $M$ are defined in Lemma \ref{l:l'M}, and $C$ is some positive constant, then $m[j_1..j_2] \geq C$ for all $j_1, j_2$ such that $w[j_1..j_2]=v$ and $r_k[j_1..j_2]\leq l'$.
\end{lemma}
{\sc Proof.} As in the previous lemma, there are no $k$-runs covering $w[i_1..i_2]$; and after we cut from $w[i_1..i_2]$ the longest prefix $w[i_1..q_1]$ covered by some left $k$-run and the longest suffix $w[q_2..i_2]$ covered by some right $k$-run, we get a factor $w[q_1+1..q_2-1]$ of measure $m[q_1+1..q_2-1]\geq 2M+2+C$ such that all $k$-runs intersecting with it are internal $k$-runs for $w[i_1..i_2]$. We illustrate the proof of Lemma \ref{l:ivp2} by Fig.~\ref{f:f1}. The ``inviolable part'' shown there is the word $w[i_1+n_1..i_1+n_2]$, whose measure $m[i_1+n_1..i_1+n_2]$ does not depend on an occurrence of $v$ to $w$.

Consider the minimal prefix $w[q_1+1..p_1]$ and the minimal suffix $w[p_2..q_2-1]$ of $w[q_1+1..q_2-1]$ such that $m[q_1+1..p_1]\geq M$ and $m[p_2..q_2-1]\geq M$. Due to Lemma \ref{l:02}, we have $m[q_1+1..p_1]\leq M+1$ and $m[p_2..q_2-1]\leq M+1$, and thus due to Lemma \ref{l:measure}, $m[p_1+1..p_2-1]\geq C$. At the same time, due to the definition of $M$, there are positions $i_1+n_1$ in $w[q_1+1..p_1]$ and $i_1+n_2$ in $w[p_2..q_2-1]$ such that $r_k(i_1+n_1)=r_k(i_1+n_2)=l'$. All the $l'$ runs contributing to $r_k(i_1+n_1)$ (resp., $r_k(i_1+n_2)$) are internal $k$-runs for $w[i_1..i_2]$ and thus do not depend on the occurrence of $v=w[i_1..i_2]$. So, in another occurrence $v=w[j_1..j_2]$ of $v$ we have
$r_k(j_1+n_1)=r_k(j_1+n_2)=l'$, and all the $k$-runs contributing to $r_k(j_1+n_1)$ (resp., $r_k(j_1+n_2)$) are internal $k$-runs for $w[j_1..j_2]$. So, there can be only internal $k$-runs for $w[j_1..j_2]$ which intersect $w[j_1+n_1..j_1+n_2]$ and thus affect its measure. So, $m[j_1..j_2]\geq m[j_1+n_1..j_1+n_2]=m[i_1+n_1..i_1+n_2]\geq m[p_1+1..p_2-1]\geq C$.

In Fig.~\ref{f:f1}, the ``inviolable part'' shown there is the word $w[i_1+n_1..i_1+n_2]$, whose measure $m[i_1+n_1..i_1+n_2]$ does not depend on an occurrence of $v$ in $w$.

\hfill $\Box$

\smallskip
The following fact is analogous to Lemma \ref{l:1+1/k}
\begin{lemma} \label{l:1+c/k}
Suppose that an infinite word $w$ satisfies the $(k,l)$-condition, and the constants $l'$ and $M$ are defined by Lemma \ref{l:l'M}. If $w[i_1..i_2]$ and $w[i_1..i_3]$ are palindromes, where $i_3>i_2$, $r_k[i_1..i_3]\leq l'$ and $m[i_1..i_2]\geq (k-1)(4l'+2M+3)$, then
\[\frac{m[i_1..i_3]}{m[i_1..i_2]}\geq 1+\frac{1}{(k-1)(4l'+2M+3)}.\]
\end{lemma}
{\sc Proof.} By Lemma \ref{l:period}, the word $w[i_1..i_3]$ is $(i_3-i_2)$-periodic. Denote its suffix $w[i_2+1..i_3]$ by $v$ and define $k'$ so that $w[i_1..i_2]=v'v^{k'-2}$ (and thus $w[i_1..i_3]=v'v^{k'-1}$) for some suffix $v'$ of $v$ not equal to $v$ (here $v'$ can be empty) and for some $k' \geq 2$.
Then $k'\leq k$ since otherwise $w[i_1..i_3]$ would have been covered by some $k$-run which is impossible due to Lemma \ref{l:covering}.

Suppose that $m[i_1..i_2]= (k-1)(4l'+2M+2+C)$ for some constant $C$. Here by the assertion we have $C\geq 1$; note that $C$ can be non-integer. Due to Lemma \ref{l:measure}, the measure of $w[i_1..i_2]$ is the sum of measures of its prefix $v'$ and of $k'-2$ occurrences of $v$. By the pigeon-hole principle, it means that some of these occurrences of $v$ (or of $v'$) have measure at least $4l'+2M+2+\lceil C \rceil$. It means by Lemma \ref{l:ivp2} that the measure $m[i_2+1..i_3]$ of another occurrence $w[i_2+1..i_3]$ of $v$ is at least $C$. So, due to Lemma \ref{l:measure} again, we have
\[\frac{m[i_1..i_3]}{m[i_1..i_2]}=1+\frac{m[i_2+1..i_3]}{m[i_1..i_2]}\geq 1+\frac{C}{(k-1)(4l'+2M+2+C)}.\]
The right hand side of this inequality is a growing function of $C$ for $C\geq 1$, so the minimal value of $C=1$ gives its minimum, and we have
\[\frac{m[i_1..i_3]}{m[i_1..i_2]}\geq 1+\frac{1}{(k-1)(4l'+2M+3)}. \phantom{aaaaa} \Box\]

\smallskip

Now we proceed to the second part of the proof, where we are going
to use Lemma \ref{l:1+c/k} to prove that there should be a factor
with palindromic length greater than $P$. The sketch of the
remaining part of the proof is the following. We assign to each
factor $w[i..j]$ its code $C[i..j]$, which is a word on a binary
alphabet. We code each factorization $\mathbb P$ of a prefix of $w[i..j]$ into palindromes
by a code $C^*([i..j],\mathbb P)$, which is obtained from $C[i..j]$ by
inserting symbols to positions between the palindromes. After that, taking a
factor $w[i..j]$ of big enough measure (the required measure depends
on $P$ and several parameters of the word), we obtain
that the number of possible codes of factorizations of prefixes of
$w[i..j]$ into at most $P$ palindromes is less than the length of
the code $C[i..j]$. We deduce from that the existence of a prefix of $w[i..j]$ not decomposable to $P$ palindromes. 

Consider a factor $w[i..j]$ with $r_k[i..j]\leq l'$ of big enough measure (``long'' factor). The required measure $N$ is determined below by \eqref{N}. 
Define the {\it
code} $C[i..j]$ as a word on the alphabet $\{1,2\}$ obtained from
the word $m(i)m(i+1)\cdots m(j)$ by erasing the symbols equal to
0. The length of the code $C[i..j]$ is denoted by $c[i..j]$;
clearly,
\begin{equation}
m[i..j]/2 \leq c[i..j] \leq m[i..j], \label {mcm}
\end{equation}
since in fact $m[i..j]$ is the sum of the symbols of $C[i..j]$, and each of them is equal to 1 or 2.

If $i$ and $j$ are fixed, consider the word $m(i)m(i+1)\cdots m(j)\in \{0,1,2\}^*$ and denote by $n_h$ the position giving the $h$th non-zero symbol in it, so that $m(n_h)\in \{1,2\}$ for all $h=1,\ldots,c$ and $m(n)=0$ for all other $n\in \{i,i+1,\ldots,j\}\backslash \{n_1,n_2,\ldots,n_c\}$. We also define $n_0=i-1$. Due to this definition, we have  $C[i..j]=a_1\cdots a_c$ with $a_h \in \{1,2\}$ and $m(n_h)=a_h$ for all $h=1,\ldots,c$.

\begin{example}
 {\rm 
Consider the prefix $w[1..25]=0101110101111111110101110$ of the Sierpinski word $w$, considered as a word satisfying the $(3,1)$-condition. Its code is $C[1..25]=111111$, and the positions $n_1,\ldots, n_6$ are equal to 4,6,10,18,22,24: these are exactly first and last positions of 3-runs in $w$, that is, of ``long'' powers of $1$. We also fix $n_0=0$.
}
\end{example}

Consider a factorization into palindromes of a prefix $w[i..i']$ of $w[i..j]$: $w[i..i']=w[i_0+1..i_1]w[i_1+1..i_2]\cdots w[i_{p-1}+1..i_p]$, where $i=i_0+1$, $i_p=i'$, and each word $w[i_d+1..i_{d+1}]$ is a palindrome. Define the {\it partition to palindromes} $\mathbb P$ as the sequence $\mathbb P=\{i_0,i_1,\ldots,i_p\}$, and the {\it code} of this partition in $w[i..j]$ as the word $C^*([i..j],i_1,\ldots,i_p)=C^*([i..j],\mathbb P)$ on the alphabet $\{1,2,*\}$ obtained from the code $C[i..j]$ by adding a star before the symbol $a_h$ for each $d$ such that $n_{h-1}\leq i_d< n_h$. Note that the number of stars in the code of a partition to $P$ palindromes is exactly $P+1$.

\begin{example}
 {\rm
Continuing the previous example, consider the following partition to palindromes of a prefix of length 22 of the Sierpinski word $w$, considered in its turn as a prefix of $w[1..25]$:
\[(010)(11)(10101111111110101)110.\]
We see that $i_0=0$, $i_1=3$, $i_2=5$, $i_3=22$, $\mathbb P=\{0,3,5,22\}$, and thus the code of this partition is
 $C([1..25],3,5,22)=C([1..25],\mathbb P)=**1*1111*1$.
}
\end{example}

The stars correspond to the boundaries of palindromes, and the symbols 1 and 2 correspond to beginnings and ends of upper $k$-runs in $w[i..j]$. Note that the code of a partition of a prefix of $w[i..j]$ always begins with a star.

We are going to prove that if $c[i..j]$ is large enough, then there is a pair of consecutive symbols from $\{1,2\}$ in $C[i..j]$ such that no partition $\mathbb P$ of a prefix of $w[i..j]$ to at most $P$ palindromes has a star between them in $C^*([i..j],\mathbb P)$. In particular it means that there is a prefix of $w[i..j]$ with palindromic length greater than $P$, and thus we prove Theorem \ref{t:klc}.


\begin{lemma}\label{l:d1d2}
 Suppose that for the word $w[i..j]$ considered above we have $c[i..j]=N$ and $i'$, $i\leq i'\leq j$ is a position such that $i'= n_h-1$ for some $h$, $1\leq h\leq N$.
 Then the number of values of $h'$ such that there exists $i''$ with $w[i'+1..i'']$ being a palindrome and $n_{h'-1}\leq i''<n_{h'}$, is bounded by $H=D_1+1+\log_{D_2}(2N/D_1)$, where $D_1=(k-1)(4l'+2M+3)$ and $D_2=1+\frac{1}{D_1}$.
\end{lemma}
{\sc Proof.} The fact that $c[i..j]=N$ means in particular that
$m[i..j]\leq 2N$ due to \eqref{mcm}. We shall estimate the number
of possible values of $m[i'+1..i'']$ not exceeding $2N$, where
$w[i'+1..i'']$ is a palindrome, and this will give an upper bound
for the number of values of $h'$, since $n_h$  are exactly the
positions where the measure changes.

First of all, $m[i'+1..i'']$ can take at most
$(k-1)(4l'+2M+3)+1=D_1+1$ values of $h'$ from $0$ to
$(k-1)(4l'+2M+3)$. Due to Lemma \ref{l:1+c/k}, the value of $h'$
numbered $D_1+2$ must be equal at least to $D_1 D_2$, where
$D_2=1+\frac{1}{(k-1)(4l'+2M+3)}=1+\frac{1}{D_1}$; and the value
of $h'$ number $D_1+n+1$ is equal at least to $D_1 D_2^n$. Even
for the maximal $n$ we should have $D_1 D_2^n \leq 2N$, so that $n
\leq \log_{D_2} (2N/D_1)$, and the total possible number of
measures of palindromes is bounded by $D_1+1+\log_{D_2}(2N/D_1).$
\hfill $\Box$

\begin{lemma}\label{l:bound}
 Suppose that for the word $w[i..j]$ considered above we have $c[i..j]=N$
 and $i'$, $i\leq i'\leq j$ is a position such that $n_{h-1}\leq i'<  n_h$ for some $h$, $1\leq h\leq N$.
 Then the number of values of $h'$, such that there exists $i''$ with $w[i'+1..i'']$ being a palindrome and $n_{h'-1}\leq i''<n_{h'}$, is bounded by $(M+4l')(\log_{D_2}(2N/D_1))+D_3$,
 where, as above, $D_1=(k-1)(4l'+2M+3)$, $D_2=1+\frac{1}{D_1}$, and $D_3=(M+4l')^2+M+4l'+D_1+1$.
\end{lemma}

{\sc Proof.} As in the previous lemma, we shall estimate the
number of possible values of $m[i'+1..i'']$ not exceeding $2N$,
where $w[i'+1..i'']$ is a palindrome. As above, $m[i'+1..i'']$ can
take at most $(k-1)(4l'+2M+3)+1=D_1+1$ values from 0 to
$(k-1)(4l'+2M+3)=D_1$.

Suppose now that $m[i'+1..i'']>D_1$, and consider the prefix
$w[i'+1..n_h-1]$ of $w[i'+1..i'']$; note that by the definition of
the sequence $\{n_h\}$, its measure $m[i'+1..n_h-1]=0$, and the
measure of $w[n_h..i'']$ is equal to that of $w[i'+1..i'']$.

Suppose first that $n_h-i'-2\geq i''-n_h$, that is, that
$w[i'+1..n_h-1]=t$ contains the center of the palindrome
$w[i'+1..i'']=u$. It means that $u=tv=\tilde{v} t' v$ and
$t=\tilde{v} t'$ for some words $v$ and $t'$. We see that the
measure of the occurrence of $\tilde{v}$ starting from $i'+1$ is
equal to 0, and thus, due to Lemma \ref{l:ivp1}, we have
$m[i'+1..i'']= m[n_h..i''] < M+4l'$, a contradiction to the fact
that $m[i'+1..i'']>D_1$.

Now suppose that $n_h-i'-2< i''-n_h$, that is, $w[i'+1..n_h-1]=t$
is shorter than a half of $w[i'+1..i'']=u$. It means that $u=t v
\tilde{t}$ for some palindrome $v$ starting from the position
$n_h$. Let us denote the measure of $v$ by $m$. The measure of the prefix
occurrence of $t$ is here equal to 0; so, due to Lemma
\ref{l:ivp1}, the measure of the suffix occurrence of $\tilde{t}$
is at most $M+4l'-1$. So, we have $m\leq m[i'+1..i'']\leq
m+M+4l'-1$, so, for each possible value of $m\geq D_1-(M+4l')$,
the measure $m[i'+1..i'']$ takes at most $M+4l'$ different values.
Adding possible values from 0 to $D_1$, we see that the total
number of values is bounded by $D_1+1+(M+4l')(H-D_1+M+4l')$, 
where
$H$ is the bound for the number of measures of palindromes
starting from $n_h$ obtained in Lemma \ref{l:d1d2}, so that
$H=D_1+1+\log_{D_2}(2N/D_1)$. Simplifying the expression, we
obtain that
$D_1+1+(M+4l')(H-D_1+M+4l')$=$(M+4l')\log_{D_2}(2N/D_1)+(M+4l')(M+4l'+1)+D_1+1$=$(M+4l')\log_{D_2}(2N/D_1)+D_3.
$
 \hfill $\Box$

\smallskip
{\sc Proof of Theorem \ref{t:klc}.}  Fix a positive integer $P$
and let $N$ be a positive integer satisfying
\begin{equation}[(M+4l')\log_{D_2}(2N/D_1)+D_3]^P<N, \label{N} \end{equation} where the
constants $D_1$--$D_3$ are defined in Lemma \ref{l:bound}.
Consider a factor $w[i..j]$ of $w$ with $r_k[i..j]\leq l'$ and
$c[i..j]
=N$; such  a factor $w[i..j]$ always exists by the
definitions of $l'$ and of $c[i..j]$. By the previous lemma, the
number of possible codes of decompositions of prefixes of
$w[i..j]$ to $P$ palindromes $w[i_d+1..i_{d+1}]$,
$d=0,\ldots,P-1$,
is at
most $[(M+4l')\log_{D_2}(2N/D_1)+D_3]^P$, and hence 
less than
$N$.
In particular it means that the position of the last star in the
code can take less than $N$ different values; but the length
$c[i..j]$ of the code of $w[i..j]$ is
 $N$. So, there
are two consecutive symbols in $c[i..j]$, corresponding to the positions
$n_h$ and $n_{h+1}$, such that the last star in the code of
the decomposition to $P$ palindromes can never appear between
them, that is, that no word $w[i..q]$, where $n_h\leq q <
n_{h+1}$, can ever be decomposed into $P$ palindromes. \hfill $\Box$

\section{Discussion} \label{s:disc}
Even if we prove that the palindromic length of factors of any aperiodic word is unbounded, some ultimately periodic words, for example, $w=(110100)^{\omega}$, contain factors having arbitrarily large
palindromic lengths. So, unlike for example the result by Mignosi, Restivo and Salemi on repetitions and periodicity \cite{mrs}, the conjectured property will not give a characterization of aperiodic words.

We prove the following property of ultimately periodic words with a uniform bound on the palindromic length of its factors:

\begin{proposition}\label{up}
Let $P$ be an integer, $w$ an ultimately periodic word such that $|u|_{\rm pal}\leq P$ for each factor $u$ of $w.$
Then $w$ has a tail $w'$ of the form $w'=v (p_1 p_2)^{\omega}$, where $p_1$ and $p_2$ are palindromes.
\end{proposition}

{\sc Proof.} Let $w'$ be a tail of $w$ having period $t.$ Consider a
factor $u$ of $w'$ with $|u|> tP$. Then $u$ can be factored as $u=v_1v_2\ldots v_m$ with $m\leq P$ and each $v_i$ a non-empty palindrome.  At least one of the palindromes $v_i$ in this factorization
has length greater than $t$. So, $t$ is a period of this long palindrome. Now the
proposition follows from the well-known fact that a period of a
palindrome has the form  $p_1 p_2$, where $p_1$ and $p_2$ are
palindromes. \hfill $\Box$

\smallskip

Proposition \ref{up} implies that if the answer to Question \ref{c1} is ``no'',
then an infinite word $w$ having a uniform bound on the palindromic length of its factors is
ultimately periodic, and moreover its period has the form $p_1
p_2$, where $p_1$ and $p_2$ are palindromes.

As we have mentioned above in Remark \ref{r1}, the $(k,l)$-conditions for some $k$ and $l$ seems to be fulfiled in particular for all morphic words, so, it is hardly probable that an example of an aperiodic word with bounded palindromic length of factors will be a morphic word. At the same time, there exist Sturmian words which do not satisfy any $(k,l)$-condition: these are exactly Sturmian words whose elements of the directive sequence are unbounded (see Chapter 2 of \cite{lothaire} for the definitions). So, this paper does not contain a proof of unbounded palindromic length of factors which would be valid for all Sturmian words.

It is also clear that 

\begin{lemma}\label{l:bin}
If there exists an aperiodic word with bounded palindromic length of factors, then there exists a binary one.
\end{lemma}
{\sc Proof.} Consider an aperiodic word $w$ on an alphabet $A=\{a_1,\ldots,a_q\}$ and for each $i \in 1,\ldots,q$ define a coding $c_i$ by $c_i(a_i)=0$, $c_i(b)=1$ for all other symbols $b \in A$. Consider the infinite words $c_1(w)$, $\ldots$, $c_q(w)$. At least one of them, let us denote it by $c(w)$, is aperiodic (since otherwise $w$ would be periodic with the period equal to the least common multiple of periods of $c_1(w)$, $\ldots$, $c_q(w)$). At the same time, for each factor $u$ of $w$ we have $|u|_{\rm pal} \geq |c(u)|_{\rm pal}$. So, if $|u|_{\rm pal}\leq P$ for all factors $u$ of $w$ and for some $P$, then the same is true for all factors $c(u)$ of the binary aperiodic infinite word $c(w)$: $|c(u)|_{\rm pal} \leq P$. \hfill $\Box$

\begin{remark}
 {\rm 
The proof of Lemma \ref{l:bin} above is valid only for the case of a finite alphabet $A$. Another proof, valid also for the infinite alphabet $\{0,1,\ldots, n,\ldots\}$, was suggested by T. Hejda who uses the morphism $c$ defined by $c(i)=10^i1$. It is not difficult to prove that $|u|_{\rm pal}\leq P$ for all factors $u$ of $w$, then $|v|_{\rm pal}\leq P+4$ for all factors $v$ of $c(w)$.
}
\end{remark}

At last, we recall again that for the case of $k$-power-free words, our proof can be extended in particular to 
 privileged words instead of palindromes (see Section \ref{s:priv}). 
However, Theorem \ref{t:klc} cannot be directly extended to privileged words since in Lemma \ref{l:bound}, 
we used the properties specific for palindromes which allowed to
apply Lemma \ref{l:ivp1}.

\section*{Acknowledgement}
The authors are grateful to S. Avgustinovich, M. Bucci, G. Fici, T. Hejda and A. De
Luca for fruitful discussions.

\end{document}